\documentclass[12pt,reqno]{amsart}

\usepackage{amsmath,amssymb,hyperref}

\newtheorem{thm}{Theorem}

\newcommand\ben{\begin{enumerate}}
\newcommand\een{\end{enumerate}}

\newcommand{\twocase}[5]{#1 \begin{cases} #2 & \text{#3}\\ #4
&\text{#5} \end{cases}   }

\newcommand\be{\begin{equation}}
\newcommand\ee{\end{equation}}
\newcommand\benn{\begin{equation*}}
\newcommand\eenn{\end{equation*}}
\newcommand\bea{\begin{eqnarray}}
\newcommand\eea{\end{eqnarray}}
\newcommand\beann{\begin{eqnarray*}}
\newcommand\eeann{\end{eqnarray*}}

\newcommand{\R}{\mathbb{R}}

\newcommand{\Z}{\mathbb{Z}}

\newcommand{\cZ}{\frak{Z}}

\newcommand{\bZ}{\mathcal{Z}}

\newcommand{\ga}{\alpha}     
\newcommand{\gl}{\lambda}    
\newcommand{\gz}{\zeta}      

\newcommand{\GL}{\operatorname{GL}}
\newcommand{\PGL}{\operatorname{PGL}}
\newcommand{\PSL}{\operatorname{PSL}}

\newcommand{\bk}{\backslash}

\begin{document}

\title[Exterior Square on $\GL(n)$]{The Dirichlet Series for the Exterior Square $L$-function on $GL(n)$}
\author{Alex V. Kontorovich} \email{alexk@math.brown.edu}
\address{Department of Mathematics,
Columbia University, New York, NY, $10027$} 
\curraddr{Department of Mathematics,
Brown University, Providence, RI, $02912$}
 \date{\today}
\maketitle
\section{Introduction}

Let $f$ be a Hecke-Maass cusp form on $\PSL(n,\Z)\bk\PGL(n,\R)/O(n)$, see e.g. \cite{Goldfeld2006}. Let $A(m_1,m_2,\dots,m_{n-1})$ be the Fourier coefficients in its Jacquet-Whittaker expansion. Then the Godement-Jacquet $L$-function associated to $f$ is given as a Dirichlet series, and in terms of the Langlands-Satake parameters $\ga_{i}(p)$, by
$$
L(f,s)=\sum_{m\ge1}{A(m,1,1,\dots,1)\over m^s}=\prod_p\prod_{i=1}^n\left(1-{\ga_i(p)\over p^s}\right)^{-1}.
$$
The exterior square $L$-function is defined via the Euler product
\be\label{fwedgeDef}
L(f,s,\wedge^2)=\prod_p\prod_{1\le i <j\le n}\left(1-{\ga_i(p)\ga_j(p)\over p^s}\right)^{-1}
.
\ee
Two distinct representations of the exterior square $L$-function are known,
the first due to Jacquet and Shalika \cite{JacquetShalika1990}, and the second discovered by
Bump and Friedberg \cite{BumpFriedberg1990}.
%
%
It is our goal in this short 
note to 
present an elementary derivation of the Jacquet-Shalika construction,
expressing the 
Euler product 
in \eqref{fwedgeDef}
as a classical Dirichlet series in the Fourier coefficients $A(m_1,\dots,m_{n-1})$.\\

On $\GL(2)$, 
$$
L(f,s,\wedge^2)=\prod_p\left(1-{\ga(p)\bar\ga(p)\over p^s}\right)^{-1} =\gz(s).
$$
On $\GL(3)$, it is easy to see that the exterior square $L$-function
$$
L(f,s,\wedge^2)=L(\tilde f,s)=\sum_m{A(1,m)\over m^s}
$$
is just the dual $L$-function corresponding to the contragredient form $\tilde f$.\\

On $\GL(4)$, experts have known for some time that the exterior square $L$-function can be expressed as a zeta function times the ``Middle'' $L$-function:
$$
L(f,s,\wedge^2)=\gz(2s) \sum_m {A(1,m,1)\over m^s}.
$$
The general formula on $\GL(n)$ 
is
 as follows.\\

\begin{thm}\label{one}
For odd $n\ge3$, the Dirichlet series for the exterior square $L$-function is given by
$$
L(f,s,\wedge^2)=\sum_{m_2,m_4,\dots,m_{n-1}\ge1}{A(1,m_2,1,m_4,1,\dots,1,m_{n-1})\over (m_2\,m_4^2\,m_6^3\cdots m_{n-1}^{(n-1)/2})^s}.
$$
For even $n\ge2$, we have
$$
L(f,s,\wedge^2)=\gz\left({n\over 2}\,s\right)\sum_{m_2,m_4,\dots,m_{n-2}\ge1}{A(1,m_2,1,m_4,1,\dots,1,m_{n-2},1)\over (m_2\,m_4^2\,m_6^3\cdots m_{n-2}^{(n-2)/2})^s}.
$$
\end{thm}


\section{Proof of Theorem \ref{one}}

As is well-known from the work of Shintani and Casselman-Shalika, the 
Fourier coefficient $A(p^{k_{1}},...,p^{k_{n-1}})$ 
is the Schur function
\be\label{AS}
A(p^{k_{1}},\dots,p^{k_{n-1}})
=
S_\lambda(\ga)
.
\ee 
Here 
$$
\ga=(\ga_{1}(p),\dots,\ga_{n}(p))
$$ 
are the Langlands-Satake parameters,
and
$$
\lambda=(\lambda_{1},\dots,\lambda_n),
\qquad
\text{ where }
\qquad
\lambda_j=\sum_{i>j} k_i
.
$$
\

Recall
the following identity \cite[(3.3)]{BumpFriedberg1990}: 
\bea\nonumber
\sum_{k_1,k_2,\dots,k_{n-1}\ge0}
S_{\gl}(\ga)\
X^{k_1+k_3+k_5+\cdots} \ 
Y^{k_2+k_3+2k_4+2k_5+\cdots} \qquad\qquad& &\\
\label{eq1}
=L_0\prod_i(1-\ga_i(p)X)^{-1}\prod_{i<j}(1-\ga_i(p)\ga_j(p)Y)^{-1},
\eea
where
$$
L_0=\twocase{}{1-Y^{n/2}}{if $n$ is even;}{1-XY^{(n-1)/2}}{if $n$ is odd.}
$$
Setting $X=p^{-s}$, $Y=p^{-w}$, using \eqref{AS}, and taking the product over all primes $p$ of both sides of \eqref{eq1} gives
$$
\cZ(s,w):=\sum_{r_1,r_2,\dots,r_{n-1}\ge1}{A(r_1,r_2,\dots,r_{n-1})\over r_1^s\, r_2^w\, r_3^{s+w}\, r_4^{2w}\, r_5^{s+2w}\cdots} 
=Z(s,w)L(f,s)L(f,w,\wedge^2),
$$
where
$$
Z(s,w)=\twocase{}{1/\gz(\frac n2 w)}{if $n$ is even;}{1/\gz(s+{n-1\over2}w)}{if $n$ is odd.}
$$

\

On the other hand, we have the following Hecke relations 
\cite[Theorem 9.3.11]{Goldfeld2006}. For $n$ odd,
\beann
A(m,1,1,\dots,1)A(1,m_2,1,m_4,1,\dots,1,m_{n-1}) \ \ \ \ \ \ \ \ \ \ \  \ \ & & \\
=\sum_{{c_2\,c_4\,c_6\cdots c_{n-1}\,c_n=m}\atop{c_2|m_2, c_4|m_4, c_6|m_6,\dots,c_{n-1}|m_{n-1}}} A(c_n,{m_2\over c_2},c_2,{m_4\over c_4},c_4,\dots,{m_{n-1}\over c_{n-1}})
.
\eeann
For $n$ even, we have
\beann
& & A(m,1,1,\dots,1)A(1,m_2,1,m_4,1,\dots,1,m_{n-2},1) \\ 
&=&\sum_{{c_2\,c_4\,c_6\cdots c_{n-2}\,c_n=m}\atop{c_2|m_2, c_4|m_4, c_6|m_6,\dots,c_{n-2}|m_{n-2}}} A(c_n,{m_2\over c_2},c_2,{m_4\over c_4},c_4,\dots,{m_{n-2}\over c_{n-2}},c_{n-2})
.
\eeann
In either case, dividing both sides by $m^s$ and $(m_2\,m_4^2\,m_6^3\cdots)^w$ and summing gives
\beann
\bZ(s,w)&:= & \left(\sum_{m\ge1}{A(m,1,1,\dots,1)\over m^s} \right)\left(\sum_{m_2,m_4,\dots\ge1}{A(1,m_2,1,m_4,1,\dots)\over (m_2\,m_4^2\,m_6^3\cdots)^w} \right) \\
&= & \sum_{m,m_2,m_4,\dots\ge1}\sum_{{c_2\,c_4\cdots=m}\atop{c_2|m_2,c_4|m_4,\cdots}}{A(c_n,{m_2\over c_2},c_2,{m_4\over c_4},c_4,\dots)\over m^s(m_2\,m_4^2\,m_6^3\cdots)^w }.
\eeann
Interchange the orders of summation and write $m_i=m_i'c_i$:
\beann
\bZ(s,w)
&= & 
\sum_{c_2,c_4,\cdots,c_n\ge1}
\sum_{{m_2',m_4',m_6'\dots\ge1}\atop{m=c_2\,c_4\cdots c_n}}{A(c_n,m_2',c_2,m_4',c_4,\dots)\over (c_2\,c_4\cdots c_n)^s(m_2'\,c_2\,m_4'^2\,c_4^2\,m_6'^3\,c_6^3\cdots)^w }.
\eeann
Rename $r_1=c_n$, $r_2=m_2'$, $r_3=c_2,\dots$, with 
$$
r_{n-1}=\twocase{}{m_{n-1}'}{if $n$ is odd;}{c_{n-2}}{if $n$ is even.}
$$ 
If $n$ is even, then 
$$
\bZ(s,w)=\cZ(s,w)
$$ 
and dividing both sides by $L(f,s)$ proves the theorem in this case. \\

If $n$ is odd, then we have an additional sum over $r_n=c_{n-1}$, which does not appear inside the Fourier coefficients. Thus
$$
\bZ(s,w)=\gz(s+{n-1\over2}w) \cZ(s,w) = L(f,s)L(f,w,\wedge^2).
$$
Again, dividing both sides by $L(f,s)$ gives the desired result.
This completes the proof. $\qquad\qquad\qquad\qquad\qquad\qquad\qquad\qquad\qquad\qquad\Box$\\

The referee has kindly pointed out to us 
that the argument above is equivalent to
the fact 
\cite[page 238
(11.9;2)]{Littlewood1940}
that 
$$
\prod_{i<j}(1-\ga_{i} \ga_{j})^{-1}=\sum S_\lambda(\ga)
,
$$ 
where the summation runs over partitions $\lambda$ whose conjugate partition is even.\\

\subsection*{Acknowledgements}
The author wishes to express his gratitude to Dorian Goldfeld,
Meera  Thillainatesan,
Sol Friedberg,
and the referee for many comments and corrections to an earlier draft.
 \\

\bibliographystyle{alpha}

\bibliography{../../AKbibliog}

\end{document}